\documentclass{amsart}%without draft
\usepackage{enumerate}
\usepackage{hyperref}

\newcommand{\norm}[1]{\left\Vert#1\right\Vert}
\newcommand{\abs}[1]{\left\vert#1\right\vert}
\newcommand{\set}[1]{\left\{#1\right\}}
\newcommand{\scp}[1]{\langle{#1}\rangle}
\newcommand{\R}{R}
\newcommand{\eps}{\varepsilon}
\newcommand{\al}{\alpha}
\newcommand{\be}{\beta}

\newcommand{\si}{\sigma}

\newcommand{\wick}{\mathop{\diamond}}

\newcommand{\cB}{\mathcal{B}}
\newcommand{\cH}{\mathcal{H}}

\newcommand{\ws}{\widetilde{\sigma}}

\newcommand{\om}{\omega}
\newcommand{\Om}{\Omega}
\newcommand{\ex}[1]{\mathsf{E}\left[\,#1\,\right]}
\newcommand{\ind}[1]{\pmb 1_{#1}}
\newcommand{\wj}{\widetilde{J}}
\newcommand{\bj}{\overline{J}}
\newcommand{\wz}{\widetilde{Z}}
\newcommand{\wx}{\widetilde{X}}

\newcommand{\yd}[1][]{Y^{\delta}_{#1}}

\begin{document}
%\doi{?}
% \issn{1029-0346}
%\issnp{1045-1129} \jvol{?} \jnum{?} \jyear{?} \jmonth{?}\received{?}
%\setcounter{section}{0} \setcounter{equation}{0}
\theoremstyle{plain}
\newtheorem{theorem}{Theorem}
\newtheorem{lemma}{Lemma}
\newtheorem{proposition}{Proposition}
\theoremstyle{definition}
\newtheorem{definition}{Definition}
\theoremstyle{remark}
\newtheorem{remark}{Remark}
%\newtheorem{exam}{Example}
%\begin{frontmatter}

\author{Yu. Mishura and G. Shevchenko}
%\address{Kyiv National Taras Shevchenko University, Mathematical
%Department, Vladimirskaya, 64, 01033, Kiev, Ukraine}
%\email{myus@univ.kiev.ua}
\thanks{The work of the first author is partially supported by
NATO grant PST.CLG.980408}

%\address{Kyiv National Taras Shevchenko University, Mathematical
%Department} \ead{zhora@univ.kiev.ua.}
%\thanks{The work is partially supported by INTAS grant YSF
%03-55-2447.}
\title[Rate of convergence of Euler approximations
of SDEs driven by fBm]{The rate of convergence of Euler
approximations for solutions of stochastic differential equations
driven by fractional Brownian motion}

\maketitle

%\keywords{ Euler approximations, stochastic differential
%equations, fractional Brownian motion, fractional white noise, the
%rate of convergence. AMS Subject Classification: 60H10, 60G18}
\begin{abstract}
The paper focuses on discrete-type approximations of solutions to
non-homogeneous stochastic differential equations (SDEs) involving
fractional Brownian motion (fBm). We prove that the rate of
convergence for Euler approximations of solutions of pathwise SDEs
driven by fBm with Hurst index $H>1/2$ can be estimated by
$O(\delta^{2H-1})$ ($\delta$ is the diameter of partition). For
discrete-time approximations of Skorohod-type quasilinear equation
driven by fBm we prove that the rate of convergence is
$O(\delta^H)$.
\end{abstract}

\section{Introduction}\label{Sec:A1.1}
Many equations which arise in modeling of processes in physics,
chemistry, biology, finance, contain randomness. This randomness is
not always well modeled by the classical Gaussian white noise
(Brownian motion) because of long-range dependence, or long memory,
of the processes under consideration. In this case the appropriate
model for the randomness is fractional Brownian motion. Recall that
$B=(B_t)_{t\geq 0}$ is called fractional Brownian motion on a
complete probability space $(\Omega, \mathcal F, P)$ with Hurst
parameter $H \in (\frac{1}{2}, 1)$ if $B$ is a centered Gaussian
process with stationary increments and covariance
$R_H(t,s)=E(B_tB_s)=\frac{1}{2}(t^{2H}+s^{2H}-|t-s|^{2H})$.

Numerical solution via time discretization of SDEs driven by
Brownian motion has  long history. We refer to the
monograph~\cite{k-p}, which contains almost complete theory of
numerical solution of such SDEs with regular coefficients. The paper
\cite{ko-pr} is devoted to Euler approximations for SDEs driven by
semimartingales. Concerning numerical solution of SDEs driven by
fBm, we mention first the paper \cite{gr-a}, where equations with
modified fBm that represents a special semimartingale are studied
(recall that fBm itself is not a semimartingale). Papers
\cite{nourdin1,nourdin2} study Euler approximations for homogeneous
one-dimensional SDEs with bounded coefficients having bounded
derivatives up to third order, driven by fBm, and prove that error
of approximation is a.s.\ equivalent to $\delta^{2H-1}\xi_t$, and
the process $\xi_t$ is given explicitly. These papers also discuss
Crank--Nicholson and Milstein schemes for SDEs driven by fBm. The
methods used by the authors cannot be applied to our case, because
they require homogeneity and high regularity of the coefficients.

We consider the stochastic differential equation on $R^d$
\begin{equation} \label{A1:1}
X_t^i=X_0^i+ \sum^m_{j=1}\int^t_0 \sigma^{ij}(s,X_s)dB^j_s+ \int^t_0
b^i (s,X_s)ds, \quad i=1,...,d,\quad t\in[0,T]
\end{equation}
where the processes $B^i, i=1,...,m$ are  fractional Brownian
motions with Hurst parameter $H$, $X_0$ is a $d$-dimensional random
variable, the coefficients $\sigma^{i\,j},
b^i:\Omega\times{[0,T]}\times{R^d}\rightarrow R$ are measurable
functions.

The integral in the right-hand side of \eqref{A1:1} can be
understood in the pathwise sense defined in \cite{Z,N} or in
Wick--Skorohod sense \cite{nua1}.
%!!!!!!!!!!!!!!!!!!!!ssylki
We treat the pathwise case first. We remind that the pathwise
integral w.r.t.\ a one-dimensional fBm $B$ can be defined as
\begin{equation*}
\int_a^b f dB = \int_a^b
\big(D^\al_{a+}f\big)(s)\big(D^{1-\al}_{b-}B_{b-}\big)(s)ds,
\end{equation*}
where
\begin{equation*}
(D^\al_{a+}f\big)(s)=\frac{1}{\Gamma(1-\al)}\bigg[\frac{f(s)}{(s-a)^\al}+
\al\int_a^s\frac{f(s)-f(u)}{(s-u)^{\al+1}}du\bigg]\ind{(a,b)}(s)
\end{equation*}
and
\begin{equation*}
\big(D^{1-\al}_{b-}B_{b-}\big)(s)=\frac{e^{-i\pi\al}}{\Gamma(\al)}\bigg[\frac{B_{b-}(s)}{(b-s)^{1-\al}}+
(1-\al)\int_s^b\frac{B_{b-}(s)-B_{b-}(u)}{(u-s)^{2-\al}}du\bigg]\ind{(a,b)}(s)
\end{equation*}
are fractional derivatives of corresponding orders,
$$
B_{b-}(s)=\big(B_s-B_b\big)\ind{(a,b)}(s).
$$
The integral exists for any $\alpha\in(1-H,\nu)$ if, for example,
$f\in C^\nu (a,b)$ with $\nu+H>1$. Moreover, in this case pathwise
integral admits an estimate
\begin{equation}\label{**}
\abs{\int_a^bfdB}\le C_0(\omega)\bigg[\int_a^b
\frac{\abs{f(s)}}{(s-a)^\al}ds+\int_a^b\int_a^s\frac{\abs{f(s)-f(u)}}{(s-u)^{\alpha+1}}du\,ds\bigg],
\end{equation}
where $C_0(\om)=
C\cdot\sup_{a<s<b}\abs{D_{b-}^{1-\al}B_{b-}(s)}<\infty$ a.s.

 Denote $\sigma=(\sigma^{ij})_{d\times
m}, b=(b^i)_{d\times 1}$ and for a matrix $A=(a^{ij})_{d\times m}$,
and a vector $y=(y^i)_{d\times 1}$ denote $|A|= \sum_{i,j} |a^{ij}|,
|y|=\sum_i |y^i|$.

We suppose that the coefficients satisfy the following assumptions
\begin{itemize}
\item[(A)] $\sigma(t,x)$ is differentiable in $x$ and there exist
such $M>0, 1-H< \beta \leq 1, \frac{1}{H}-1<\kappa \leq 1$ and for
any $N>0$ there exists such $M_N>0$  that
\begin{itemize}
\item[1)] $|\sigma(t,x)-\sigma(t,y)|\leq M|x-y|, x,y \in R^d, t\in
[0,T]$; \item[2)]
$|\partial_{x_i}\sigma(t,x)-\partial_{x_i}\sigma(t,y)| \leq
M_N|x-y|^{\kappa},|x|,|y|\leq N$, $t\in[0,T]$; \item[3)]
$|\sigma(t,x)-\sigma(s,x)|+|\partial_{x_i}\sigma(t,x)-\partial_{x_i}\sigma(s,x)|
\leq M|t-s|^{\beta}, x\in R^d, t,s \in[0,T]$.
\end{itemize}
\item[(B)]
\begin{itemize} \item[1)] for any $N>0$ there exists $L_N>0$ such that $$|b(t,x)-b(t,y)| \leq
L_N |x-y|,\ |x|,|y| \leq N,\  t\in[0,T];$$ \item[2)] $|b(t,x)|\leq
L(1+|x|)$.%, $x\in\R^d$,\ $t\in[0,T]$.
\end{itemize}
\end{itemize}

As it was stated in \cite{N}, under conditions (A)--(B) the equation
\eqref{A1:1} has the unique solution $\{X_t$, $t \in [0,T]\}$, and
for a.a.\ ${\omega} \in\Omega$ this solution belongs to
$C^{H-\rho}[0,T]$ for any $0<\rho<H$. Now, let $t \in
[0,T],\delta=\frac{T}{N}, \tau_n=\frac{nT}{N}=n\delta, n=0,...,N$.
Consider discrete Euler approximations of solution of equation
\eqref{A1:1},
\begin{equation*}
\widetilde{Y}_{\tau_{n+1}}^{i,\delta}=\widetilde{Y}_{\tau_{n}}^{i,\delta}
+
b^i(\tau_n,\widetilde{Y}_{\tau_{n}}^{\delta})\delta+\sum^m_{j=1}\sigma^{ij}
(\tau_n,\widetilde{Y}_{\tau_{n}}^{\delta})\Delta B_{\tau_n}^j,\quad
\widetilde{Y}_{0}^{i,\delta}=X_0^i,
\end{equation*}
and corresponding continuous interpolations
\begin{equation}\label{contint}
{Y}_{t}^{i,\delta}=\widetilde{Y}_{\tau_{n}}^{i,\delta} +
b^i(\tau_n,\widetilde{Y}_{\tau_{n}}^{\delta})(t-\tau_n)+\sum^m_{j=1}\sigma^{ij}
(\tau_n,\widetilde{Y}_{\tau_{n}}^{\delta})(B_{t}^{j}-B_{\tau_n}^j),\quad
t\in[\tau_n,\tau_{n+1}].
\end{equation}

Continuous interpolations satisfy the equation
\begin{equation}\label{continteq}
Y_t^{i,\delta}=X_0^i+\int^t_0
b^i(t_u,Y^{\delta}_{t_u})du+\sum_{j=1}^m
\int^t_0\sigma^{i\,j}(t_u,Y^{\delta}_{t_u})dB^j_u,
\end{equation}
where $t_u = \tau_{n_u}$, $n_u=\max\{n:\tau_n\leq u\}.$

For simplicity we denote the vector of solutions as
$X_t=(X^i_t)_{i=1,...,d},$ vector of continuous approximations as
$Y_t^{\delta}=(Y_t^{\delta, i})_{i=1,...,d}$. Throughout the paper,
$C$ denotes a generic constant, whose value is not important and may
change from line to line, and we write  $C(\cdot),$ if the
dependence on some parameters is crucial.

The paper is organized as follows. Sections 2 and 3 are devoted to
equations with pathwise integral. Section 2 describes the growth and
H\"older properties of approximations $Y_t^\delta$. We use here
growth and H\"older estimates of solution of corresponding pathwise
equations from the paper \cite{N}. Section 3 contains estimates of
rate of convergence for Euler approximations of the solutions of
pathwise equations. It is well-known that in the case $H=\frac12$,
when we have SDE with It\^o integral with respect to Wiener process,
the rate of convergence is $O(\delta^{1/2})$. It is natural to
expect that in our case the rate of convergence might be
$O(\delta^H)$. Nevertheless, the rate of convergence is only of
order $O(\delta^{2H-1})$ unless the diffusion coefficient is
constant (where the rate is $O(\delta^H)$, as expected).
%A simple
%example in Section 4 demonstrates that even in the linear case the
%rate of convergence is exactly $O(\delta^{2H-1})$, which can be
%considered as fractional analogue of well-known ``Clark--Cameron
%paradox'' (see \cite{cl-ca}).
Section 4 is devoted to discrete-time approximations of solutions of
SDEs with stochastic divergence integral with respect to the fBm
(otherwise known as Skorohod integral, or fractional white noise
integral). It is shown that the rate of convergence is $O(\delta^H)$
in this case. The better rate of convergence is mainly because in
the pathwise case there is no ``It\^o compensator'' for the
integral.

%%%%%%%%%%%%%%%%%%%%%%%%%%%%%%%%%%%%%%%%%%%%%%%%%%%%%%%%%%%%%

%\sectionmark{Existence conditions}
\section{Some properties of Euler approximations for solutions of pathwise equations}\label{Sec:A1.2}
%\sectionmark{Existence conditions}
In this section we consider growth and H\"older properties of the
approximation process $\set{\yd[t],t\in[0,T]}$. We need some
additional notations. Denote
$\varphi_{u,v}:=\abs{\yd[t_u]-\yd[v]}(u-v)^{-\alpha-1}$ for
$0<v<t_u<T$, $0<\alpha<1$, $X_t^* := \sup_{0\le s\le t}\abs{X_s}$,
$Y_t^{\delta,*}:=\sup_{0\le s\le t}\abs{\yd[s]}$. Further, for any
$0<\rho<H$ there exists such $C=C(\omega, \rho)$ that for any
$0<v<u$
\begin{equation} \label{eq3} \abs{B_u-B_v}\le C(\omega,
\rho)(u-v)^{H-\rho}.
\end{equation}
 We shall use the following statement \cite[Lemma 7.6]{N}
\begin{proposition}\label{P1}
Let $0<\alpha<1$, $a,b>0$, $x:\R_+\to R_+$ be a continuous function
such that for each $t$
$$
x_t\le a+ b t^\alpha \int_0^t (t-s)^{-\alpha} s^{-\alpha} x_s ds.
$$
Then $x_t\le a c_\alpha \exp\set{d_\al t b^{1/(1-\al)}} $, where
$c_\al= 4e^2\frac{\Gamma(1-\alpha)}{1-\al}$, $d_\al =
2\big(\Gamma(1-\al)\big)^{1/(1-\al)}$, $\Gamma(\cdot)$ is Euler's
Gamma function.
\end{proposition}
We also establish technical lemma, which will be used later.
\begin{lemma}\label{alem}
There exists such $C=C_\al>0$ that for any $s\in[0,T]$, $s\neq t_s$
and $\delta\le 1$, $\alpha\in(0,1)$ it holds
$$
J:=\int_0^{t_s} (s-u)^{-\al-1}\int_u^{t_u} (v-t_v)^{-\al}dv\,du\le
C\delta^{-\al}.
$$
\end{lemma}
\begin{proof}
Evidently,
\begin{align*}
J&=\int_0^{t_s}(v-t_v)^{-\al}\int_0^v (s-u)^{-\al-1}du\,dv\le
\al^{-1} \int_0^{t_s} (v - t_v)^{-\al}(s-v)^{-\al}dv.
\end{align*}
Let $t_s = n\delta$ for some $0<n\le N$. Then
$$
\int_0^{t_s}(v - t_v)^{-\al}(s-v)^{-\al}dv =
\sum_{k=0}^{n-2}\int_{\tau_k}^{\tau_{k+1}}+
\int_{(n-1)\delta}^{(2n-1)\delta/2} +
\int_{(2n-1)\delta/2}^{n\delta}.
$$
We estimate the integrals individually:
\begin{align*}
&\int_{\tau_k}^{\tau_{k+1}} \le (s-\tau_{k+1})^{-\al}
\int_{\tau_k}^{\tau_{k+1}} (v-t_v)^{-\al}dv \le
(1-\al)^{-1}(s-\tau_{k+1})^{-\al}\delta^{1-\al},\\
&\int_{(n-1)\delta}^{(2n-1)\delta/2}\le
(\delta/2)^{-\al}\int_{(n-1)\delta}^{(2n-1)\delta/2}
(v-t_v)^{-\al}dv\le C\delta^{1-2\al},\\
& \int_{(2n-1)\delta/2}^{n\delta} \le
(\delta/2)^{-\al}\int_{(2n-1)\delta/2}^{n\delta}(s-v)^{-\al}dv\le
C\delta^{1-2\al}.
\end{align*}
Therefore
$$\aligned
J&\le C\delta^{1-2\al}+\delta^{-\al}\sum_{k=0}^{n-2}
(s-\tau_{k+1})^{-\al}\delta\le
C\delta^{1-2\al}+\delta^{-\al}\int_0^{n\delta} (s-v)^{-\al}dv\\&\le
C\delta^{1-2\al}+C\delta^{-\al}\le C\delta^{-\al}. \endaligned
$$
\end{proof}
\eject
%
%
% We start with growth and H\"{o}lder properties of the approximation $Y_t^{\delta}.$
\begin{theorem}\label{T1}
(i) Let the conditions {\rm (A)--(B)}  hold and

{\rm (C) 1)}\qquad $|\sigma(t,x)| \leq C(1+|x|)$.

Then for any $\eps>0$ and $0<\rho<H$ there exists $\delta_0>0$ and
$\Omega_{\eps,\delta_0,\rho}\subset \Omega$ such that
$P(\Omega_{\eps,\delta_0,\rho})>1-\eps$ and for any
$\om\in\Om_{\eps,\delta_0,\rho}$, $\delta<\delta_0$ one has
$\abs{\yd[t]}\le C(\omega)$, $\abs{\yd[t_s]-\yd[t_r]}\le
C(\omega)(t_s-t_r)^{H-\rho}$, $0\le r<s\le T$.

(ii) If, instead of {\rm (A), 2)} and {\rm (C)} we assume that $b$
and $\si$ are bounded functions, then $\abs{\yd[t]}\le C(\omega)$,
$\abs{\yd[s]-\yd[r]}\le C(\omega)(s-r)^{H-\rho}$, $0\le r<s\le T$.
%\end{itemize}
%\end{itemize}

In both cases $C(\omega)$ does not depend on $\delta$.
\end{theorem}

\begin{proof}
We can assume that $\delta\le 1$. It follows immediately from (A),
1) and 3) and \eqref{continteq} that for any $\alpha\in
(1-H,\beta\wedge 1/2)$
\begin{equation}
\label{eq4}
\begin{aligned}
|Y_t^{i,\delta}|&\le \abs{X_0^i} + \int_0^t
\abs{b^i(t_u,\yd[t_u])}du
+ \sum_{j=1}^m \abs{\int_0^t\si^{i\,j}(t_u,\yd[t_u])dB_u^H}\\
& \le \abs{X_0^i} + L\int_0^t\Big(1+\abs{\yd[t_u]}\Big)du +
C_0(\om)\sum_{j=1}^m \int_0^t
\abs{\si^{i\,j}(t_u,\yd[t_u])}u^{-\al}du \\
&\quad  + C_0(\om)\sum_{j=1}^m \int_0^t\int_0^r
\abs{\si^{i\,j}(t_r,\yd[t_r])-\si^{i\,j}(t_u,\yd[t_u])}(r-u)^{-\al-1}du\,dr\\
&\le \abs{X_0^i} + \Big(C_0(\om) \frac{T}{1-\al}+LT\Big) +
\big(C_0(\om)+CT^\al\big) \int_0^t\abs{\yd[t_u]}u^{-\al}du\\
&\quad +
MC_0(\om)\int_0^t\int_0^{t_r}\Big((t_r-t_u)^\be+\abs{\yd[t_r]-\yd[u]}+\abs{\yd[u]-\yd[t_u]}\Big)(r-u)^{-\al-1}
du\,dr.
\end{aligned}
\end{equation}
(We use here the equality $t_r=t_u$ for $t_r\le u<r$.) Denote
$C_1(\om):=m\big(C_0(\om)\frac{T^{1-\al}}{1-\al}+LT\big)+\abs{X_0}$,
$C_2(\om):= m(C_0(\om)+CT^\al)$. Further, note that $t_r-t_u\le
r-u+\delta$. Also, it follows from representations \eqref{contint}
 that for any $\rho\in(0,H)$
\begin{equation}
\label{eq5}
\begin{aligned}
\abs{\yd[u]-\yd[t_u]}&\le L\Big(1+\abs{\yd[t_u]}\Big)(u-t_u)+C\cdot
C(\om,\rho)\Big(1+\abs{\yd[t_u]}\Big)(u-t_u)^{H-\rho}\\
&\le C_3(\om)\Big(1+\abs{\yd[t_u]}\Big)(u-t_u)^{H-\rho},
\end{aligned}
\end{equation}
where $C_3(\om)=LT^{1-H-\rho}+C\cdot C(\om,\rho)$.

Moreover, for $\be>\al$
\begin{align*}
P_t&:=\int_0^t\int_0^{t_r}(t_r-t_u)^\beta(r-u)^{-\al-1}du\,dr\le\int_0^t\int_0^{t_r}\big((r-u)^\beta+\delta^\be\big)(r-u)^{-\al-1}du\,dr\\
&\le (\be-\al)^{-1}\int_0^t r^{\be-\al}dr +
\al^{-1}\delta^\be\int_0^t(r-t_r)^{-\al}dr,
\end{align*}
and for any $k\ge 0$ and any power $\pi>-1$
$$
\int_{\tau_k}^{\tau_{k+1}}
(r-t_r)^{\pi}dr=\int_{\tau_k}^{\tau_{k+1}}
(r-\tau_k)^{\pi}dr=C_1\delta^{\pi+1}\text{ with }C_1=(\pi+1)^{-1},
$$
whence
\begin{equation}
\label{eq6} \int_0^t (r-t_r)^{-\al}dr\le \int_0^T (r-t_r)^{-\al}dr=
C_1N\delta^{1-\al}=C_1\delta^{-\al}.
\end{equation}
Therefore
\begin{equation}
\label{eq7} P_t\le C_1T^{\be-\al+1}+\al^{-1}C_1\delta^{\be-\al}\le
C_1T^{\be-\al+1}+\al^{-1}C_1=: C_2.
\end{equation}
Estimate now $$ Q_t:=\int_0^t\int_0^{t_r}
\abs{\yd[u]-\yd[t_u]}(r-u)^{-\al-1}du\,dr,
$$
using \eqref{eq5} and \eqref{eq6}:
\begin{equation}\label{eq8}
\begin{aligned}
Q_t&\le
\big(1+Y_t^{\delta,*}\big)\int_0^t\int_0^{t_r}(u-t_u)^{H-\rho}(r-u)^{-\al-1}du\,dr\\
&\le
C_3(\om)\big(1+Y_t^{\delta,*}\big)\delta^{H-\rho}\alpha^{-1}\int_0^t(r-t_r)^{-\al}dr\le
C_4(\om)\big(1+Y_t^{\delta,*}\big)\delta^{H-\alpha-\rho},
\end{aligned}
\end{equation}
with $C_4(\om)=C_3(\om)\alpha^{-1}\cdot C_1$. Note that
$Y_t^{\delta,*}:=\sup_{0\le s\le t} \abs{\yd[s]}<\infty$ for any
$t\in[0,T]$ a.s. Substituting \eqref{eq7} and \eqref{eq8} into
\eqref{eq4}, we obtain that
\begin{equation}
\label{V}
\begin{aligned} \abs{\yd[t]}& \le
C_5(\om)+C_2(\om)\int^t_0\abs{\yd[t_u]}u^{-\al}du+
C_4(\om)\big(1+Y_t^{\delta,*}\big)\delta^{H-\al-\rho}\\&\quad+C_6(\om)\int_0^t\int_0^{t_r}\varphi_{r,u}du\,dr
\end{aligned}
\end{equation}
with $C_5(\om)=C_3(\om)+MC_0(\om)C_2$, $C_6(\om)=MC_0(\om)$. To
simplify the notations, in what follows we remove subscripts from
$C(\om)$,  writing $C(\omega)$ for all constants depending  on
$\om$.

So we can write
\begin{equation}
\label{eq9}
\begin{aligned}
Y^{\delta,*}_t&\le
C(\om)\Big(1+Y_t^{\delta,*}\delta^{H-\alpha-\rho}+\int_0^t\abs{\yd[t_u]}u^{-\alpha}du+\int_0^t\int_0^{t_r}\varphi_{r,u}du\,dr\Big).
\end{aligned}
\end{equation}
In turn, we can estimate $\int_0^{t_s}\varphi_{s,u}du$. At first,
similarly to the previous estimates,
\begin{equation}
\label{eq10}
\begin{aligned}
\abs{\yd[t_s]-\yd[u]}&\le
C(\om)\Big[\int_u^{t_s}\Big(1+\abs{\yd[t_v]}\Big)dv+\int_{u}^{t_s}\Big({1+\abs{\yd[t_v]}}\Big)(v-u)^{-\alpha}dv\\
&\quad +
\int_u^{t_s}\int_u^{t_v}\abs{\si(t_v,\yd[t_v])-\si(t_z,\yd[t_z])}(v-z)^{-\alpha-1}dz\,dv\Big]\\
&\le
C(\om)\Big[(t_s-u)^{1-\al}+\int_u^{t_s}\abs{\yd[t_v]}(v-u)^{-\al}dv
+
\delta^\be \int_u^{t_s} (v-t_v)^{-\al}dv\\
&\quad +\int_u^{t_s}\int_u^{t_v}\varphi_{v,z}dz\,dv +
\int_u^{t_s}\int_u^{t_v}\abs{\yd[z]-\yd[t_z]}(v-z)^{-\al-1}dz\,dv\Big];
\end{aligned}
\end{equation}
multiplying by $(s-u)^{-\al-1}$ and integrating over $[0,t_s]$, we
obtain that
\begin{equation}
\label{eq11} \int_0^{t_s} \varphi_{s,u}du\le C(\om)\sum_{i=1}^5
Q_s^i,
\end{equation}
where
\begin{align}\label{eq12} Q_s^1&:=\int_0^{t_s}
(t_s-u)^{1-\al}(s-u)^{-\al-1}du\le \int_0^{t_s}(s-u)^{-2\al}du\le
C;\\
%\begin{aligned}
Q_s^2&:=\int_0^{t_s}
(s-u)^{-\al-1}\int_u^{t_s}\abs{\yd[t_v]}(v-u)^{-\al}dv \label{eq13}\\
& = \int_0^{t_s} \abs{\yd[t_v]}\int_0^v
(v-u)^{-\alpha}(s-u)^{-\al-1}du\,dv\le
C_0\int_0^{t_s}\abs{\yd[t_v]}(s-v)^{-2\al}dv,\nonumber
%\end{aligned}
\end{align}
where $C_0=\int_0^\infty (1+y)^{-\al-1}y^{-\al}dy$; according to
Lemma~\ref{alem}
\begin{equation}\label{eq14}\begin{aligned}
Q_s^3& := \delta^\be \int_0^{t_s}
(s-u)^{-\al-1}\int_u^{t_s}(v-t_v)^{-\al}dv\,du\\
&%\le
%\delta^\be\int_0^{t_s}(s-u)^{-\al}du\int_0^{t_s}(v-t_v)^{-\al}dv
\le
C\delta^\be\delta^{-\al}\le C.
\end{aligned}
\end{equation}
%The last integral can be estimated as
%\begin{align*}
%&\int_0^{t_s}(s-v)^{-\al}(v-t_v)^{-\al}dv\le
%\int_0^{t_s}(t_s-v)^{-\al}(v-t_v)^{-\al}dv\\
%&\le
%\int_{t_s-\delta}^{t_s} (t_s-v)^{-\al}(v-t_s)^{-\al}dv+\int_0^{t_s-\delta} \delta^{-\al}(v-t_s)^{-\al}dv\\
%=C
%\end{align*}
Further, using estimates \eqref{eq5}, we can conclude that
\begin{equation}
\label{eq16}
\begin{aligned}
Q^4_s&:=\int_0^{t_s}(s-u)^{-\al-1}\int_u^{t_s}\int_u^{t_v}\varphi_{v,z}dz\,dv\,du\\
&\le \int_0^{t_s} \int_0^{t_v}\int_0^{z\wedge
v}\varphi_{v,z}(s-u)^{-\al-1}du\,dz\,dv\le C\int_0^{t_s}
(s-v)^{-\al}\int_0^{t_v}\varphi_{v,z}dz\,dv.
%\\
%&\quad+ C(\om)
%\int_0^{t_s}(s-v)^{-\alpha}\int_{t_v}^{v}(z-t_v)^{H-\rho}(v-z)^{-1}dz\,dv
%\Big(1+\abs{Y^{\delta,*}_{t_s}}\Big)\\
%&\le C\int_0^{t_s} (s-v)^{-\al}\int_0^{t_v}\varphi_{v,z}dz +
%C(\om)\Big(1+\abs{Y^{\delta,*}_{t_s}}\Big)\delta^{H-\al-\rho}.
\end{aligned}
\end{equation}
At last, using estimates (7) and Lemma~\ref{alem}, we can conclude
that.
\begin{equation}
\label{eq17}
\begin{aligned}
Q_s^5&:=\int_0^{t_s}(s-u)^{-\al-1}\int_u^{t_s}\int_u^{t_v}\abs{\yd[z]-\yd[t_z]}(v-z)^{-\alpha-1}dz\,dv\,du\\
&\le
C(\om)\int_0^{t_s}(s-u)^{-\al-1}\int_u^{t_s}\int_u^{t_v}(v-z)^{-\al-1}dz\,dv\,du\cdot\delta^{H-\rho}\Big(1+\abs{Y^{\delta,*}_{t_s}}\Big)\\
&\le
C(\om)\Big(1+\abs{Y^{\delta,*}_{t_s}}\Big)\delta^{H-\rho-\alpha}.
\end{aligned}
\end{equation}
Now, denote $\psi_s:=Y^{\delta,*}_s+\int_0^{t_s}\varphi_{s,u}du$.
Note that the integrals $Q_s^i$ are finite for $s=k\delta$, i.e. for
any $s\in[0,T]$, including $s=t_s$. Then it follows from \eqref{eq9}
and \eqref{eq11}--\eqref{eq17} that
$$
\psi_t\le C(\om)\Big(1+Y^{\delta,*}_t\delta^{H-\alpha-\rho} +
\int_0^t\big((t-v)^{-2\alpha}+v^{-\alpha}\big)\psi_v dv\Big).
$$
 Let $\eps>0$ be fixed. Note that all
constants $C(\om)$ are finite a.s.\ and independent of $\delta$.
Thus, we can choose $\delta_0>0$ and $\Om_{\eps,\delta_0,\rho}$ such
that $C(\omega)\delta_0^{H-\alpha-\rho}\le 1/2$ on
$\Om_{\eps,\delta_0,\rho}$ and $P(\Om_{\eps,\delta_0,\rho})>1-\eps$.
Then for any $\om\in\Om_{\eps,\delta_0,\rho}$
$$
\psi_t\le
C(\om)+\frac12\psi_t+C(\om)\int_0^t\big((t-v)^{-2\alpha}+v^{-\alpha}\big)\psi_v
dv,
$$
whence
$$
\psi_t\le C(\omega)\Big(1+t^{2\al}\int_0^t
(t-v)^{-2\al}v^{-2\al}\psi_v dv\Big),
$$
and it follows immediately from the last equation and
Proposition~\ref{P1} that $\psi_t\le C(\om)$ whence, in particular,
$\abs{\yd[t]}\le C(\om)$, $t\in[0,T]$, and $\int_0^{t_s} \varphi_u
du\le C(\om)$. Moreover, from \eqref{eq10} with $u=t_r$, $r\le s$,
taking into account that $\int_{t_r}^{t_s}(v-t_v)^{-\al}dv\le
\delta^{-\alpha}(t_s-t_r)$, we obtain the estimate
\begin{align*}
\abs{\yd[t_s]-\yd[t_r]}&\le
C(\om)\Big((t_s-t_r)^{1-\al}+\delta^{\be-\al}(t_s-t_r)+(t_s-t_r)\\&\quad+\delta^{H-\rho}\int_{t_r}^{t_s}(v-t_v)^{-\al}dv\Big)\le
C(\om)(t_s-t_r)^{1-\al},
\end{align*}
and the statement (i) is proved.
%!!!!!!!!!!!!!!!!!!!!!!!!!!!!!!!!abs
(ii) Let $\abs{b(t,x)}\le b$, $\abs{\si(t,x)}\le\si$. Then it is
very easy to see that the estimate \eqref{V} will take a form
$$
\abs{\yd[t]}\le
C(\om)\Big(1+\int_0^t\int_0^{t_r}\varphi_{r,u}du\,dr\Big),
$$
\eqref{eq10} will perform to
\begin{align*}
\abs{\yd[t_s]-\yd[u]}&\le
C(\om)\Big((t_s-u)^{1-\al}+(\delta^\be+\delta^{H-\rho})\int_u^{t_s}(v-t_v)^{-\al}dv\\
&\quad +\int_u^{t_s}\int_u^{t_v} \varphi_{v,z}dz\,dv \Big)
\end{align*}
and instead of \eqref{eq11}--\eqref{eq17} we obtain $$
\int_0^{t_s}\varphi_{s,u}du\le
C(\om)\Big(1+\int_0^{t_s}(s-v)^{-\al}\int_0^{t_v}\varphi_{v,z}dz\,dv\Big),
$$
whence the proof easily follows.
\end{proof}

%!!!!!!!!!!!!!!!!!!!!!!!!!!!!!!!!!!!!t_s
%\begin{remark}
%We established the statement for psi..... Besov
%\end{remark}

%%%%%%%%%%%%%%%%%%%%%%%%%%%%%%%%%%%%%%%%%%%%%%%%%%%%%%%%%%%%%

%\sectionmark{Existence conditions}
\section{The estimates of rate of convergence for Euler approximations of the solutions
of pathwise equations}\label{Sec:A1.3}
%\sectionmark{Existence conditions}

Now we establish the estimates of the rate of convergence of our
approximations \eqref{continteq} for the solution of the equation
\eqref{A1:1} with pathwise integral w.r.t.\ fBm. We establish even
more, namely, the estimate of convergence rate for the norm of the
difference $X_t-\yd[t]$ in some Besov space, similarly to the result
of Theorem~1. Denote
$$
\Delta_{u,s}(X,\yd):=\abs{X_s-\yd[s]-X_u+\yd[u]}
$$
and assume for technical simplicity that $L_N=L$, $M_N=M$ in (A) and
(B).
\begin{theorem}\label{t2}
Let the conditions {\rm(A)--(C)} hold and also
\begin{itemize}
\item[\rm(D)]
\begin{itemize} \item[\rm 1)] H\"{o}lder continuity of the coefficient $b$ in time:
$|b(t,x)-b(s,x)|\leq C|t-s|^{\gamma}$, $C>0$, $2H-1<\gamma\leq 1$;
\item[\rm 2)] the exponent $\beta$ from {\rm (A) 3)} satisfies
$\beta>H$.
 \end{itemize} \end{itemize}
 Then:

 (i) for any $\eps>0$ and any $\rho>0$ sufficiently small there exists $\delta_0>0$ and
 $\Om_{\eps,\delta_0,\rho}$ such that
 $P(\Om_{\eps,\delta_0,\rho})>1-\eps$ and for any
 $\om\in\Om_{\eps,\delta_0,\rho}$, $\delta<\delta_0$
$$
U_\delta:=\sup_{0\le s\le
T}\Big(\abs{X_s-\yd[s]}+\int_0^{t_s}\abs{\Delta_{u,s}(X,\yd)}
(s-u)^{-\al-1}du\Big)\le C(\om)\cdot \delta^{2H-1-\rho},
$$
where $C(\om)$ does not depend on $\delta$ and $\eps$ (but depends
on $\rho$);

(ii) if, in addition, the coefficients $b$ and $\si$ are bounded,
then for any $\rho\in(0,2H-1)$  there exists $C(\om)<\infty$ a.s.\
such that $U_\delta\le C(\om)\delta^{2H-1-\rho}$, $C(\om)$ does not
depend on $\delta$.
\end{theorem}
\begin{proof}
(i) Denote $Z_t^\delta:=\sup_{0\le s\le t}\abs{X_s-\yd[s]}$. Then
\begin{equation}\label{eq-8}
\begin{split}
&Z^\delta_t:=\sup\limits_{0\leq s\leq
t}|X_s-Y^\delta_s|\leq\sup\limits_{0\leq s\leq
t}\int_0^s|b(u,X_u)-b(t_u,Y^\delta_{t_u})|du\\&+ \sup\limits_{0\leq
s\leq
t}\sum\limits_{i,j=1}^m|\int_0^s(\sigma^{ij}(u,X_u)-\sigma^{ij}(t_u,Y^\delta_{t_u}))
dB_u^i|\leq \int_0^t|b(u,X_u)-b(u,\yd[u])|du
\\
&+\int_0^t|b(u,\yd[u])-b(t_u,\yd[u])|du+
\int_0^t|b(t_u,\yd[u])-b(t_u,
Y^\delta_{t_u})|du\\&+\sup\limits_{0\leq s\leq t}
\sum\limits_{i,j=1}^m|\int_0^s(\sigma^{i\,j}(u,X_u)-\sigma^{i\,j}(u,
\yd[u]))dB_u^i|
\\&+\sup\limits_{0\leq s\leq t}
\sum\limits_{i,j=1}^m|\int_0^s(\sigma^{i\,j}(u,\yd[u])-\sigma^{ij}(t_u,
\yd[u]))dB_u^i| \\&+\sup\limits_{0\leq s\leq
t}\sum\limits_{i,j=1}^m|\int_0^s(\sigma^{i\,j}(t_u,\yd[u])-
\sigma^{i\,j}(t_u,Y^\delta_{t_u}))dB_u^i|=: \sum\limits_{k=1}^6I_k.
\end{split}
\end{equation}
Now we estimate separately all these terms. Evidently,
\begin{equation}
\label{eq-9} I_1\le L\int_0^t Z_u^\delta du.
\end{equation}
Condition (D) 1) implies that for $\delta\le1$
\begin{equation}
\label{eq-10} I_2\le C\int_0^t \abs{u-t_u}^\gamma du\le
C\delta^\gamma\le C\delta^{2H-1}.
\end{equation}
 As it follow from Theorem~2.2, for any
$\eps>0$ and any $\rho\in(0,H)$ there exists $\delta_0>0$ and
$\Om_{\eps,\delta_0,\rho}\subset\Om$  such that
$P(\Om_{\eps,\delta_0,\rho})>1-\eps$ and $C(\om)$ independent of
$\eps$ and $\delta$ such that for any $\om\in
\Om_{\eps,\delta_0,\rho}$ it holds $\abs{\yd[t]-\yd[s]}\le
C(\om)\abs{t-s}^{H-\rho}$. In what follows we assume that
$\delta<\delta_0<1$. Therefore
\begin{equation}
\label{eq-11} I_3\le L\cdot C(\om)\delta^{H-\rho}\cdot t\le
C(\om)\delta^{H-\rho}, \om\in \Om_{\eps,\delta_0,\rho}.
\end{equation}
Now we go on with $I_4$. It follows from \eqref{**} that for
$1-H<\al<1/2$
\begin{equation}\label{eq-12}
\begin{aligned}
&I_4  \le C(\om)\sum_{i,j=1}^m\bigg[\int_0^t
\abs{\si^{i\,j}(u,X_u)-\si^{i\,j}(u,\yd[t_u])}u^{-\al}du\\
&+\int_0^t \int_0^r
\abs{\si^{i\,j}(r,X_r)-\si^{i\,j}(u,X_u)-\si^{i\,j}(r,\yd[r])+
\si^{i\,j}(u,\yd[u])}\\&\qquad\qquad\times(r-u)^{-\al-1}du\,dr
\bigg]=:I_7+I_8.
\end{aligned}
\end{equation}
Evidently,
\begin{equation}\label{eq-13}
I_7\le C(\om)\int_0^t Z_u^\delta u^{-\al} du.
\end{equation}
According to \cite[Lemma 7.1]{N}, under condition (A)
\begin{equation}\label{eq-14}
\begin{aligned}
&\abs{\si(t_1,x_1)-\si(t_2,x_2)-\si(t_1,x_3)+\si(t_2,x_4)} \le
M\abs{x_1 - x_2 - x_3 + x_4}\\&\qquad + M\abs{x_1 -
x_3}\Big(\abs{t_2-t_1}^\be + \abs{x_1-x_2}^\kappa +
\abs{x_3-x_4}^\kappa\Big).
\end{aligned}
\end{equation}
Therefore, $I_8\le\sum_{k=9}^{12} I_k$, where
\begin{align*}
I_9&= C(\om) \int_0^t \int_0^r \abs{X_r
-\yd[r]}(r-u)^{\be-\al-1}du\,dr,\\
I_{10}&=C(\om)\int_0^t \int_0^r
\abs{X_r-\yd[r]}\abs{X_r-X_u}^\kappa(r-u)^{-\al-1}du\,dr,\\
I_{11}&=C(\om)\int_0^t \int_0^r
\abs{X_r-\yd[r]}\abs{\yd[r]-\yd[u]}^\kappa(r-u)^{-\al-1}du\,dr,\\
I_{12}&=C(\om)\int_0^t \int_0^r
\Delta_{u,r}(X,\yd)(r-u)^{-\al-1}du\,dr.
\end{align*}
Taking into account that $\be>H>\al$, we obtain that
\begin{equation}\label{eq-15}
I_9\le C(\om)\int_0^t Z^\delta_u du.
\end{equation}
As it follows from \cite[Theorem 2.1]{N}, under assumptions (A) and
(B) for any $0<\rho<H$ there exists such constant $C(\om)$ that
\begin{equation}\label{eq-16}
\sup_{0\le t\le T} \abs{X_t}\le C(\om),\quad \sup_{0\le s\le t\le
T}\abs{X_t-X_s}\le C(\om)\abs{t-s}^{H-\rho}.
\end{equation}
Moreover, we can choose $\rho>0$ and $\al>1-H$ such that
$\kappa(H-\rho)>\al$ and $H-\rho>2H-1$, because $\kappa H > 1-H$. In
this case
\begin{equation}
\label{eq-17} I_{10}\le C(\om) \int_0^t Z^\delta_r \int_0^r
(r-u)^{\kappa(H-\rho)-\al-1}du\,dr\le C(\om) \int_0^T Z_r^\delta dr.
\end{equation}
It follows from Theorem 2.2 that on $\Om_{\eps,\delta_0,\rho}$ the
same estimate holds for $I_{11}$.

Now estimate $I_5$.
\begin{align*}
&I_5 \le C(\om) \int_0^t \abs{\si(u,\yd[u])-\si(t_u,\yd[u])}u^{-\al}
du \\
&+ C(\om)\int_0^t \int_0^r
\abs{\si(r,\yd[r])-\si(t_r,\yd[r])-\si(u,\yd[u])+\si(t_u,\yd[u])}(r-u)^{-\al-1}du\,dr\\
&=: I_{13}+I_{14}.
\end{align*}
Obviously, %\begin{subequations}
\begin{align}
I_{13}&\le C(\om)\delta^\be, \label{eq17.1}\\ I_{14}&\le
C(\om)\Big(\int_0^t\int_0^{t_r}+\int_0^t\int_{t_r}^r\Big)\abs{...}(r-u)^{-\al-1}du\,dr\nonumber\\
&\le C(\om)\int_0^t\int_0^{t_r}\delta^\be (r-u)^{-\al-1}du\,dr+
\int_0^t\int_{t_r}^r\big((r-u)^\be+(r-u)^{H-\rho}\big) du\,dr
\nonumber\\&\le
C(\om)\big(\delta^{\be-\al}+\delta^{H-\rho-\al}\big)\label{eq17.2}.
\end{align}
Similarly,
\begin{equation}
\label{eq17.3}
\begin{aligned}
&I_6\le C(\om) \int_0^t
\abs{\si(t_u,\yd[u])-\si(t_u,\yd[t_u])}u^{-\al}du\\
& + C(\om) \int_0^t  \int_0^r
\abs{\si(t_r,\yd[r])-\si(t_r,\yd[t_r])-\si(t_u,\yd[u])+\si(t_u,\yd[t_u])}\\
&\qquad\qquad\times(r-u)^{-\al-1}du\,dr=:I_{15}+I_{16}.
\end{aligned}
\end{equation}
Here
\begin{align}
I_{15}&\le C(\om)\int_0^t \delta^{H-\rho} u^{-\al}du\le
C(\om)\delta^{H-\rho};\label{eq17.3}\\
I_{16}&\le C(\om)\int_0^t \int_0^r \delta^{H-\rho}
(r-u)^{-\al-1}du\,dr\le C(\om)\delta^{H-\rho-\al}.\label{eq17.4}
\end{align}
Substituting \eqref{eq-9}--\eqref{eq17.4} into \eqref{eq-8}, we
obtain that on $\Om_{\eps,\delta_0,\rho}$
\begin{equation}
\label{eq-18} Z_t^\delta \le C(\om) \Big(\int_0^tZ_r^\delta r^{-\al}
dr + \delta^{H-\rho-\al}+\delta^{H-\rho}+\int_0^t\theta_r dr\Big),
\end{equation}
where $\theta_r=\int_0^r \Delta_{r,u}(X,\yd)(r-u)^{-\al-1}du$.
Recall that $H-\rho>2H-1$, therefore
$$
Z_t^\delta \le C(\om) \Big(\int_0^t\big(Z_r^\delta
r^{-\al}+\theta_r\big)dr+ \delta^{2H-1-\rho} \Big).
$$
We now estimate $\theta_r$. Evidently, for $t>u$
\begin{align*}
\Delta_{t,u}(X,\yd)&\le\int_u^t\abs{b(s,X_s)-b(t_s,\yd[t_s])}ds\\&\qquad
+ \sum_{i,j=1}^m \abs{\int_u^t
\big(\si^{i\,j}(s,X_s)-\si^{i\,j}(t_s,\yd[t_s])\big)dB_s^i}.
\end{align*}
Therefore, using inequality \eqref{**}, we obtain that $\theta_t\le
\sum_{k=1}^9 J_k$, where
\begin{align*}
J_1&= \int_0^t\int_u^t
\abs{b(s,X_s)-b(s,\yd[s])}ds(t-u)^{-\al-1}du,\\
J_2&= \int_0^t\int_u^t
\abs{b(s,\yd[s])-b(t_s,\yd[s])}ds(t-u)^{-\al-1}du,\\
J_3&=\int_0^t\int_u^t
\abs{b(t_s,\yd[s])-b(t_s,\yd[t_s])}ds(t-u)^{-\al-1}du,\\
J_4&=C(\om)\int_0^t\int_u^t
\abs{\si(s,X_s)-\si(s,\yd[s])}(s-u)^{-\al}ds(t-u)^{-\al-1}du,\\
J_5&=C(\om)\int_0^t\int_u^t
\abs{\si(s,\yd[s])-\si(t_s,\yd[s])}(s-u)^{-\al}ds(t-u)^{-\al-1}du,\\
J_6&=C(\om)\int_0^t\int_u^t
\abs{\si(t_s,\yd[s])-\si(t_s,\yd[t_s])}(s-u)^{-\al}ds(t-u)^{-\al-1}du,\\
%J_7&=C(\om)\int_0^t\int_u^t
%\abs{b(t_s,\yd[s])-b(t_s,\yd[t_s])}ds(t-u)^{-\al-1}du,\\
J_7&=C(\om)\int_0^t\int_u^t\int_u^r
\abs{\si(r,X_r)-\si(r,\yd[r])-\si(v,X_v)+\si(v,\yd[v])}\\&\qquad\qquad\times(r-v)^{-\al-1}dv\,dr(t-u)^{-\al-1}du,\\
J_8&=C(\om)\int_0^t\int_u^t\int_u^r
\abs{\si(r,\yd[r])-\si(t_r,\yd[r])-\si(v,\yd[v])+\si(t_v,\yd[v])}\\&\qquad\qquad\times(r-v)^{-\al-1}dv\,dr(t-u)^{-\al-1}du,\\
J_{9}&=C(\om)\int_0^t\int_u^t\int_u^r
\abs{\si(t_r,\yd[r])-\si(t_r,\yd[t_r])-\si(t_v,\yd[v])+\si(t_v,\yd[t_v])}\\&\qquad\qquad\times(r-v)^{-\al-1}dv\,dr(t-u)^{-\al-1}du.
\end{align*}
It is clear that
\begin{align*}
J_1\le C\int_0^t Z_s^\delta \int_0^s(t-u)^{-\al-1}du\,ds,\quad
J_2\le C\delta^\gamma,\quad J_3\le C(\om)\delta^{H-\rho}.
\end{align*}
Further,
$$
J_4\le C\int_0^tZ_s^\delta\int_0^s(s-u)^{-\al}(t-u)^{-\al-1}du\,ds.
$$
As we noted before,  the inner integral
$\int_0^s(s-u)^{-\al}(t-u)^{-\al-1}du\le C_0(t-s)^{-2\al}$,
$C_0=\int_0^\infty (1+y)^{-\al-1}y^{-\al}dy$. Therefore
$$
J_4\le C\int_0^t(t-s)^{-2\al} Z_s^\delta ds.
$$
Similarly to $J_2$, $J_5\le C(\om)\delta^\gamma$, and similarly to
$J_3$, $J_6\le C(\om)\le C(\om)\delta^{H-\rho}$. Estimating  $J_7$,
$J_8$ and $J_9$ is, of course, a bit more complicated, but not
dramatically. Obviously,
\begin{align*}
J_8&\le C(\om) \delta^\be \int_0^t\int_u^t\int_u^r
(r-v)^{-\al-1}dv\,dr(t-u)^{-\al-1}du\\&= C(\om)\delta^\be\int_0^t
(t-u)^{-2\al}du\le C(\om)\delta^\be;
\end{align*}
similarly $J_9\le C(\om)\delta^{H-\rho}$. Now we apply to $J_7$ the
inequality \eqref{eq-14} and obtain the following estimate of the
integrand:
\begin{equation}
\label{eq-19}
\begin{aligned}
& \abs{\si(r,X_r)-\si(r,\yd[r])-\si(v,X_v)+\si(v,\yd[v])} \le
M\Big[\Delta_{r,v} (X,\yd) \\& +
\abs{X_r-\yd[r]}(r-v)^\be+\abs{X_r-\yd[r]}\abs{X_r-X_v}^\kappa +
\abs{X_r-\yd[r]}\abs{\yd[r]-\yd[v]}^\kappa\Big].
\end{aligned}
\end{equation}
According to this, we write $J_7\le \sum_{k=10}^{13}J_k$, where, in
turn,
\begin{align*}
J_{10}& = C(\om) \int_0^t\int_u^t\int_u^r \Delta_{r,v} (X,\yd)
(r-v)^{-\al-1}dv\,dr(t-u)^{-\al-1}du \\
& = C(\om) \int_0^t\int_0^r\int_0^v (t-u)^{-\al-1}du\Delta_{r,v}
(X,\yd) (r-v)^{-\al-1}dr\,dv\\
&\le C(\om)\int_0^t(t-r)^{-\al}\theta_r dr;\\
J_{11}&=C(\om)\int_0^t\int_u^t\int_u^r \abs{X_r-\yd[r]}
(r-v)^{\be-\al-1}dv\,dr(t-u)^{-\al-1}du\\
&\le C(\om)\int_0^t Z_r^\delta\int_0^r(t-u)^{-\al-1}\Big(\int_u^r
(r-v)^{\be-\al-1} dv \Big) du\,dr\\
&\le C(\om)\int_0^t (t-r)^{-\al}Z_r dr,\\
J_{12}& = C(\om)\int_0^t\int_u^t\int_u^r
\abs{X_r-\yd[r]}\abs{X_r-X_v}^\kappa
(r-v)^{-\al-1}dv\,dr(t-u)^{-\al-1}du\\
&\le C(\om)\int_0^t\int_0^r\int_u^r Z_r^\delta
(r-v)^{\kappa(H-\rho)-\al-1}dv(t-u)^{-\al-1}du\,dr\\
&\le C(\om)\int_0^t Z_r^\delta (t-r)^{-\al}dr,
\end{align*}
and $J_{13}\le C(\om)\int_0^t Z_r^\delta (t-r)^{-\al}dr$ is obtained
the same way. Summing up these estimates, we obtain that
$$
J_{7}\le C(\om)\int_0^t (t-r)^{-\al}
\big(Z_r^\delta+\theta_r\big)dr,
$$
whence \begin{equation} \label{eq-20} \theta_t\le C(\om)
\Big(\int_0^t(t-r)^{-2\al}\big(Z_r^\delta+\theta_r\big)dr+
\delta^{H-\rho}+\delta^\gamma\Big).
\end{equation}
Coupling  together \eqref{eq-18} and \eqref{eq-20}, and taking into
account that $H-\rho>2H-1$, $\gamma> 2H-1$, we obtain
\begin{equation}\label{eq-21}
\begin{aligned}
Z_t^\delta+\theta_t &\le C(\om)\Big(\delta^{2H-1} +\int_0^t
\big((t-r)^{-2\al}+r^{-\al}\big)\big(Z_r^\delta+\theta_r\big)dr\Big)\\
&\le C(\om)\Big(\delta^{2H-1} +t^{2\al}\int_0^t
(t-r)^{-2\al}r^{-2\al}\big(Z_r^\delta+\theta_r\big)dr\Big)
\end{aligned}
\end{equation}
The proof now follows immediately from \eqref{eq-21} and
Proposition~2.1.

The statement (ii) is obvious.
\end{proof}
\begin{remark}%------------------------------------!
In \cite{nourdin2} it is proved that for an equation with
homogeneous regular coefficients the error
$\abs{X_t-\yd[t]}\delta^{1-2H}$ almost surely converges to some
stochastic process $\xi_t$, which means that the estimate of the
rate of convergence in Theorem~\ref{t2} is sharp.
\end{remark}
%------------------------------------------------------

\section{Approximation of quasilinear Skorohod-type equations}
Now we proceed to the problem of numerical solution of Skorohod-type
equation driven by fractional white noise. From now on, we assume
that our probability space is the white noise space $(\Omega,
\mathcal{F},P)=(S'(\R),\cB(S'(\R)),\mu)$, $\diamond$ is the Wick
product, $B^0_t=\langle \om, \ind{[0,t]}\rangle$ is Brownian motion,
$W^0=\dot{B}^0$ is the white noise (see \cite{chto-to} for
definitions). Next, in order to introduce an fBm with Hurst
parameter $H>1/2$ on this space, we define for $f:[0,T]\to \R$ the
fractional integral operator
$$
M f(x)=K \int_x^T (s-x)^{H-3/2}f(s)\,ds, $$ where
$$
\gathered K=\bigl(\sin(\pi
H)\Gamma(2H+1)\bigr)^{1/2}\bigl(K_1^2+K_2^2\bigr)^{-1/2},
\\
K_1=\pi\Bigl(2\cos((3/4- H/2)\pi)\Gamma(3/2-H)\Bigr)^{-1},
\\
K_2=\pi\Bigl(2\sin((3/4- H/2)\pi)\Gamma(3/2-H)\Bigl)^{-1},
\endgathered
$$
and set $M_{t}(x)=M \ind{[0,t]} (x).$ We also define for
$f,g:[0,T]\to \R $ the scalar product and the norm
\begin{gather*}
\scp{f,g}_{\mathcal{H}}=H(2H-1) \int_0^T\int_0^T
f(t)g(s)\abs{t-s}^{2H-2}dt\,ds, \quad
\norm{f}_{\mathcal{H}}^2=\scp{f,f}_{\mathcal{H}}.
\end{gather*} %,\quad \alpha=H(2H-1).
 The process
$$
B_t=\langle{M_t,\om}\rangle,\quad t\in[0,T] $$ is the fBm with Hurst
parameter $H$. Let also $W=\dot B$ be the fractional white noise.
Detailed description of the white noise theory can be found in
\cite{evdfh}, \cite{huoks}.

Consider quasilinear Skorohod-type equation driven by fractional
white noise
\begin{equation}
\label{main SDE}
X(t)=X_0+\int_0^t b(s,X(s),\om)\,ds+\int_0^t \sigma(s) X(s)\wick
W(s) \,ds
\end{equation}
with non-random initial condition $X_0$. Suppose that coefficients
$b$ and $\sigma$ satisfy the following conditions:
%\begin{subequations}
%\label{cond1}
\begin{itemize}
\item[(E)]
\begin{itemize}
\item [1)] The linear growth condition and Lipschitz condition on $b$:
\begin{equation*}%\label{cond1a}
\begin{gathered}
\abs{b(t,x,\om)}\le C(1+\abs{x}),\\
\abs{b(t,x,\om)-b(t,y,\om)}\le C\abs{x-y};
\end{gathered}
\end{equation*}
\item[2)] ``Smoothness'' of $b$ w.r.t.\ $\om$: for any $t\in[0,T]$ and for
$h\in L^1(\R)$
\begin{equation*}
%\label{cond1b}
\abs{b(t,x,\om+h)-b(t,x,\om)}\le C(1+\abs{x})\int_\R\abs{h(s)}ds.
\end{equation*}
% b(t,x,\om)$ стохастично диференційовний і його похідна має лінійний рост
%по $x$:
%\begin{equation}\label{cond1b}
%\abs{D_s b(t,x,\om)}\le C(1+\abs{x});
%\end{equation}
\item[3)]H\"older continuity of $b$ w.r.t.\ $t$ or order $H$ with constant that
grows linearly in $x$:
\begin{equation*}%\label{cond1c}
\abs{b(t,x,\om)-b(s,x,\om)}\le C(1+\abs{x})\abs{t-s}^{H};
\end{equation*}
\item[4)] H\"older continuity of $\sigma$ w.r.t.\ $t$ or order $H$:
\begin{equation*}%\label{cond1d}
\abs{\sigma(t)-\sigma(s)}\le C\abs{t-s}^{H}.
\end{equation*}
\end{itemize}
\end{itemize}
\begin{remark}
As in previous sections, we denote by $C$ any constant which may
depend on coefficients of the equation, on initial condition $X_0$
and on the time horizon $T$, but is independent of anything else
(and we write $C(v)$ to emphasize the dependence on $v$).
\end{remark}
\begin{remark}
The condition~{\rm (E) 2)} is true if, for example, the coefficient
$b$ has stochastic derivative growing at most linearly in $x$. It is
obviously true if $b$ is non-random.
\end{remark}
 Define for $t\in[0,T]$ $\sigma_t(s)=\sigma(s)\ind{[0,t]}(s)$ and
denote
$$
J_{\sigma}(t)=\exp^{\diamond}\set{-\int_0^t
\sigma(s)dB_s}=\exp\set{-\int_{\R} M \sigma_t (s) dB^0(s) - \frac12
\norm{\sigma_t}^2_{\cH} ds}
$$
the fractional Wick exponent. It follows from \cite[Theorem
2]{misquasilin} that under assumptions (E) equation~\eqref{main SDE}
has the unique solution that belongs to all $L^p$ and can be
represented in the form
\begin{equation*}
X(t)=J_{-\sigma}(t)\wick Z(t),
\end{equation*}
where the process $Z(t)$ solves (ordinary) differential equation
\begin{equation}\label{varofcon}
Z(t)=X_0+\int_0^t
J_\sigma(s)b(s,J^{-1}_{\sigma}(s)Z(s),\om+M\sigma_s)\,ds.
\end{equation}
This gives the following idea of constructing time-discrete
approximations of the solution of \eqref{main SDE}. Take the
uniform partitioning $\set{\tau_n=n\delta,\ n=1,\dots,N}$ of the
segment $[0,T]$ and define first the approximations of $Z$ in a
recursive way:
\begin{equation}
\label{vocappr}
\begin{aligned}
&\wz(0)=X_0,\\
&\wz(\tau_{n+1})=\wz(\tau_n)+
\wj(\tau_n)
b(\tau_n,\wj^{-1}(\tau_n) \wz(\tau_n),\om+M\ws_{n})\delta,
\end{aligned}
\end{equation}
where
\begin{equation*}
\begin{aligned}
& \wj(t)=\exp\set{-\int_0^t \ws(s)dB_s-\frac12 \norm{\ws\ind{[0,t]}}^2_{\cH}},\\
& \ws(s)= \sigma(t_s),\ \ws_n=\ws\ind{[0,\tau_n]}.%\\
%& m_n = M \ws_n
\end{aligned}
\end{equation*}
Note that both
$\norm{\ws_n}_{\cH}$ and $M\ws_n$ are easily computable as
finite sums of elementary integrals% (but see
%Remark~\ref{simul} below)
.  Further, we interpolate continuously by
\begin{equation}\label{appr4z}
\wz(t)=X_0+\int_0^t \wj(t_s) b(t_s,
\wj^{-1}(t_s)\wz(t_s),\om+M\ws_{{n_s}})\, ds,
\end{equation} where $n_s = \max\{n:\tau_n\le s\}$,
%!!!!!!!!!!!!!!!!!!!!!!!!!!!!!!!!!!!!n_s======
and set
\begin{equation}\label{appr4x}
\wx(t) = T_{-M\ws\ind{[0,t]}}\wj^{-1}(t)\wz(t),
\end{equation}
where for $h\in S'(\R)$ $T_h$ is the shift operator, $T_h F(\om)=
F(\om+h)$.

\begin{lemma}
\label{l3} Under the assumption {\rm (E) 1)} the following estimate
is true
\begin{equation*}
\abs{e^{\alpha_1}b(t,e^{-\alpha_1}x,\om)-e^{\alpha_2}b(t,e^{-\alpha_2}x,\om)}
\le C(1+e^{\alpha_1}+e^{\alpha_2}+\abs{x})\abs{\alpha_1-\alpha_2}.
\end{equation*}
\end{lemma}
\begin{proof}
Write
\begin{align*}
&\abs{e^{\alpha_1}b(t,e^{-\alpha_1}x,\om)-e^{\alpha_2}b(t,e^{-\alpha_2}x,\om)}\\
&\ \le \abs{e^{\alpha_1}b(t,e^{-\alpha_1}x,\om)-e^{\alpha_1}b(t,e^{-\alpha_2}x,\om)}
+\abs{e^{\alpha_1}b(t,e^{-\alpha_2}x,\om)-e^{\alpha_2}b(t,e^{-\alpha_2}x,\om)}
\end{align*}
and apply {\rm (E) 1)}.
\end{proof}
\begin{lemma}
\label{l1}
Let $\xi_1$ and $\xi_2$ be jointly Gaussian variables. Then for
$q\ge 1$
$$
\ex{\abs{e^{\xi_1}-e^{\xi_2}}^{2q}}\le
C(L,q)\left(\ex{(\xi_1-\xi_2)^2}\right)^q,
$$
where $L=\max\set{\ex{\xi_1^2},\ex{\xi_2^2}}$.
\end{lemma}
\begin{proof}
By Lagrange theorem, Cauchy--Schwartz inequality and Gaussian
property,
\begin{align*}
\ex{\abs{e^{\xi_1}-e^{\xi_2}}^{2q}}&\le
\left(\ex{e^{4q\xi_1}+e^{4q\xi_2}}
\ex{\abs{\xi_1-\xi_2}^{4q}}\right)^{1/2}\\&\le
C(L)C(q)\left(\ex{(\xi_1-\xi_2)^2}\right)^q,
\end{align*}
as required.
\end{proof}

Our first result is about convergence of $\wz$ to $Z$.
\begin{theorem}\label{thz1}
Under conditions {\rm (E)} for any $p\ge 1$ the following estimate
holds:
\begin{equation}
\ex{\abs{Z(t)-\wz(t)}^{2p}}\le C(p)\delta^{2pH}.
\end{equation}
\end{theorem}
\begin{proof}
Firstly, we remind that $Z(t)$ belongs to all $L^q$ and
$\ex{\abs{Z(t)}^q}\le C(q)$. Therefore equation \eqref{varofcon}
together with the condition {\rm (E) 2)} gives
$\ex{\abs{Z(t)-Z(s)}^q}\le C(q)\abs{t-s}^q$. Equation
\eqref{vocappr} and the condition {\rm (E) 1)} allow to write
$$\abs{\wz(\tau_{n+1})}\le (1+C\delta)\abs{\wz(\tau_n)}+C\delta \wj(\tau_n)\le
e^{C\delta}\abs{\wz(\tau_n)}+C\delta \wj(\tau_n).$$ This gives an
estimate
$$
\abs{\wz(\tau_n)}\le C\sum_{k=0}^{N-1} \wj (\tau_k)\delta.
$$
Then for any $q\ge 1$ by the Jensen inequality,
$$
\abs{\wz(\tau_n)}^q \le C(q)\sum_{k=0}^{N-1} \wj^q(\tau_k)\delta,
$$
Taking expectations, we get $$\ex{\abs{\wz(\tau_n)}^q}\le
C(q)\sum_{k=0}^{N-1}\ex{\wj^q(\tau_k)}\delta.$$ Using that each
$\wj$ is exponent of Gaussian variable and $\sigma$ is bounded on
$[0,T]$, we obtain
$$
\ex{\abs{\wz(\tau_n)}^q}\le
C(q)\sum_{k=0}^{N-1}\delta = C(q).
$$
This through \eqref{appr4z} and {\rm (E) 1)}  implies
$\ex{\abs{\wz(t)}^q}\le C(q)$.
% and in turn
%$\ex{\abs{\wz(t)-\wz(s)}^q}\le C_q \abs{t-s}^q$. Consequently, it
%is enough to prove
%\begin{equation*}
%%\end{equation*}

Now write
\begin{equation*}
\abs{Z(t)-\wz(t)} \le I_1+I_2+I_3+I_4+I_5,
\end{equation*}
where
\begin{align*}
& I_1 = \left\lvert \int_0^t \wj(t_s)\big(b(t_s,\wj^{-1}(t_s)
Z(t_s),\om+M\ws_{n_s})\right.\\
&\hspace{15em}\left.-\vphantom{\int_0^t}b(t_s,\wj^{-1}(t_s)
\wz(t_s),\om+M\ws_{n_s})\big)\,ds\right\rvert,\\
& I_2 = \left\lvert\int_0^t \big(\wj(t_s)
b(t_s,\wj^{-1}(t_s)Z(t_s),\om+M\ws_{n_s})\right.\\
&\hspace{15em} \left.\vphantom{\int_0^t} - J_{\sigma}(s)
b(t_s,J_{\sigma}^{-1}(s)Z(t_s),\om+M\ws_{n_s})\big)\,
ds \right\rvert,\\
& I_3 = \left\lvert\int_0^t J_{\sigma}(s)\big(b(s,J_{\sigma}^{-1}(s)
Z(t_s),\om+M\ws_{n_s})\right.\\
&\hspace{15em} \left.\vphantom{\int_0^t} -b(t_s,J_{\sigma}^{-1}(s)
Z(t_s),\om+M\ws_{n_s})\big)\,ds\right\rvert,\\
& I_4 = \abs{\int_0^t J_{\sigma}(s)\big(b(s,J_{\sigma}^{-1}(s)
Z(t_s),\om+M\ws_{n_s})-b(s,J_{\sigma}^{-1}(s)
Z(t_s),\om+M\sigma_s)\big)\,ds},\\
& I_5 = \abs{\int_0^t J_{\sigma}(s)\big(b(s,J_{\sigma}^{-1}(s)
Z(s),\om+M\sigma_s)-b(s,J_{\sigma}^{-1}(s)
Z(t_s),\om+M\sigma_s)\big)\,ds}.
%& A_5 = \abs{\int_0^t J_{\sigma}(s)\big(b(s,J_\sigma^{-1}(s)
%\wz(s),\om+\sigma_t)-b(s,\wj^{-1}(t_s)
%\wz(s),\om+\sigma_t)\big)\,ds};\\
%& A_5 = \abs{\int_0^t J_{\sigma}(s)\big(b(s,J_\sigma^{-1}(s)
%Z(s),\om+\sigma_t)-b(s,J_{\sigma}^{-1}(s)
%\wz(s),\om+\sigma_t)\big)\,ds}.
\end{align*}
We first estimate using Lemma~\ref{l3}%(and
%boundedness of $\sigma$ which follows from the last)
\begin{align*}
I_2&\le
C\int_0^t\big(1+J_\sigma(s)+\wj(t_s)+\abs{Z(t_s)}\big)\left(\abs{\int_0^s
\big(\sigma(u)-\ws(u)\big)\,dB_u}\right.\\
&\left.\vphantom{\int_0^s}\hspace{9em}+\abs{\sigma(t_s)\big(B_s-B(t_s)\big)}+
\frac12\abs{\norm{\sigma_s}^2_{\cH}-\norm{\ws_{n_s}}^2_{\cH}} \right)ds\\
&\le C\int_0^t\big(1+J_\sigma(s)+\wj(t_s)+\abs{Z(t_s)}\big)\\
&\hspace{9em}\cdot\left(\abs{\int_0^s
\big(\sigma(u)-\ws(u)\big)\,dB_u}+\abs{B_s-B_{t_s}}+\delta^{H}\right)ds,%!!!!!!!!!!!!!!!!!!!
\end{align*}
where the inequality
$\abs{\norm{\sigma_s}^2_{\cH}-\norm{\ws_{n_s}}^2_{\cH}}<C\delta^H$
is due to E 4) and boundedness of $\sigma$ on $[0,T]$. Applying
Cauchy--Schwartz inequality, we arrive to
\begin{align*}
I_2\le &
C\left(\int_0^T\big(1+J_{\sigma}^2(s)+\wj^2(t_s)+Z^2(t_s)\big)ds
\right)^{1/2}\\
& \qquad \cdot\left( \int_0^T \biggl(\Big(\int_0^s
\big(\sigma(u)-\ws(u)\big)\,dB_u\Big)^2+ (B_t-B_{t_s})^2+\delta^{2H}
\biggr)ds  \right)^{1/2}.
\end{align*}
Further, from {\rm (E) 3)}
\begin{equation*}
I_3\le C \int_0^T (J_\sigma(s)+\abs{Z(s)})\,ds \delta^{H},
\end{equation*}
from {\rm (E) 2)}
$$
I_3\le C \int_0^T (J_\sigma(s)+\abs{Z(s)})\,ds \delta^{H} .
$$
Condition {\rm (E) 1)} allows to estimate
\begin{align*}
&I_1 \le C\int_0^t \abs{Z(t_s)-\wz(t_s)}\,ds,\\
&I_5\le C\int_0^t\abs{Z(s)-Z(t_s)}\,ds.
\end{align*}
Summing up these estimates yields
\begin{align*}
&\abs{Z(t)-\wz(t)}\le
C\left(\int_0^T\big(1+J_{\sigma}^2(s)+\wj^2(t_s)+Z^2(t_s)\big)ds
\right)^{1/2}\\
&\quad  \cdot\left(\delta^{2H}+ \int_0^T \biggl(\Big(\int_0^s
\big(\sigma(u)-\ws(u)\big)\,dB_u\Big)^2+ (B_t-B_{t_s})^2
\biggr)ds  \right)^{1/2}\\
&\ +C\int_0^T
\abs{Z(t_s)-\wz(t_s)}ds+C\int_0^t\abs{Z(s)-Z(t_s)}\,ds.
\end{align*}
Then, using (discrete) Gronwall inequality, we get
\begin{align*}
&\abs{Z(t)-\wz(t)}\le
C\left(\int_0^T\big(1+J_{\sigma}^2(s)+\wj^2(t_s)+Z^2(t_s)\big)ds
\right)^{1/2}\\
&\quad \cdot\left(\delta^{2H}+ \int_0^T \biggl(\Big(\int_0^s
\big(\sigma(u)-\ws(u)\big)\,dB_u\Big)^2+ (B_t-B_{t_s})^2
\biggr)ds  \right)^{1/2}\\
&\ +C\int_0^t\abs{Z(s)-Z(t_s)}\,ds.
\end{align*}
Then we raise this to the $2p$th power and use Jensen's inequality.
The last term will be bounded by $C(p)\delta^{2p}$, in the first one
we apply Cauchy--Schwartz inequality for expectations, Jensen's
inequality and use uniform boundedness of moments for $Z$,
$J_{\sigma}$ and $\wj$ (for $J_{\sigma}$ and $\wj$ it follows from
the fact that the both are exponents of some Gaussian variables with
bonded variance) to get
\begin{align*}
\ex{\abs{Z(t)-\wz(t)}^{2p}}&\le C(p)
\biggl(
\delta^{2pH}+
 \Bigl(
 \mathsf{E}{\Big[\,\Bigl\lvert\int_0^T \big(\sigma(u)-\ws(u)\big)\,dB_u\Bigr\rvert^{4p}\,\Big]}
 \Bigr)^{1/2}
\\ &\hspace{10em}
 + \vphantom{\int_0^T}
 \left(
 \ex{\abs{B_t-B_{t_s}}^{4p}}
 \right)^{1/2}
\biggr).
\end{align*}
Using again that $\ex{\abs{\cdot}^{4p}}=C(p)(\ex{(\cdot)^2})^{2p}$
for Gaussian variables, we get
\begin{align*}
\ex{\abs{Z(t)-\wz(t)}^{2p}}&\le C(p)
\biggl(
\delta^{2pH}+
 \Bigl(
 \mathsf{E}\Big[\,{\Big\lvert\int_0^T \big(\sigma(u)-\ws(u)\big)\,dB_u\Big\rvert^2}\,\Big]
 \Bigr)^{p}
\\
& \hspace{10em}+\vphantom{\int_0^T}
 \left(
 \ex{\abs{B_t-B_{t_s}}^{2}}
 \right)^{p}
\biggr)\\
& \le C(p) \bigl(\delta^{2pH}+\norm{\sigma-\ws}^{2p}_{\cH}\bigr)\le C(p)\delta^{2pH},
\end{align*}
the last is due to {\rm (E) 4)}. This is the desired result.
\end{proof}
Now we are ready to state the main result of this section.
\begin{theorem}\label{th52}
Under conditions {\rm (E)} approximations $\wx$ defined by
\eqref{appr4x} converge to the solution $X$ of \eqref{main SDE} in
the mean-square sense, and moreover
$$
\ex{(X(t)-\wx(t))^2}\le C\delta^{2H}.
$$
\end{theorem}
\begin{proof}
Estimate first for $h\in L^1(\R)$
\begin{align*}
&{T_{h} Z(t)-Z(t)}\le A_1+A_2+A_3\\
&A_1= \int_0^t {T_h J_\sigma(s)}
\Big\lvert b(s,(T_h J^{-1}_\sigma)T_h Z(s),\om+h+M\sigma_s)\\
& \hspace{15em}-
 b(s,(T_h J^{-1}_\sigma)Z(s),\om+h+M\sigma_s)\Big\rvert\,ds,\\
&A_2=\int_0^t  {T_h J_{\sigma}(s)}
\Big\lvert b(s,(T_hJ_\sigma^{-1})Z(s),\om+h+M\sigma_s)\\
& \hspace{15em} -b(t,(T_h J_\sigma^{-1}(s))
Z(s),\om+M\sigma_s)\Big\rvert\,ds,\\
&A_3=\int_0^t \Big\lvert T_h
J_\sigma(s)b(t,(T_h J_\sigma^{-1}(s)) Z(s),\om+M\sigma_s)\\
& \hspace{15em}-J_\sigma(s)b(t,J_\sigma^{-1}(s)
Z(s),\om+M\sigma_s)\Big\rvert\,ds.%\\
%\le & C\int_0^t\abs{T_hZ(s)-Z(s)}\,ds + C\int_0^T(1+\abs{Z(s)})\abs{h(s)}ds
%+C \int_0^t (1+\abs{Z(s)})\left(\right)
\end{align*}
The condition {\rm (E) 1)} gives $A_1\le
C\int_0^t\abs{T_hZ(s)-Z(s)}\,ds$, the condition {\rm (E) 2)} gives
$$
A_2\le C\int_0^T(1+\abs{Z(s)})\,ds \int_{\R}\abs{h(s)}ds
$$
and Lemma~\ref{l3} with boundedness of $\sigma$ yields
\begin{align*}
A_3&\le C\int_0^T(1+J_{\sigma}(s)+T_hJ(\sigma)+\abs{Z(s)})\,ds
\abs{\int_{\R}M \sigma(s)\,h(s)\,ds} . \\
&\le C\int_0^T (1+J_{\sigma}(s)+T_hJ(\sigma)+\abs{Z(s)})\,ds
\int_{\R}\abs{h(s)}ds.
\end{align*}
Applying Gronwall lemma, we get
\begin{align*}
\abs{T_hZ(t)-Z(t)}\le & C\int_0^T
(1+J_{\sigma}(s)+T_hJ(\sigma)+\abs{Z(s)})\,ds \int_{\R}\abs{h(s)}ds.
\end{align*}
Raising this inequality to the $2p$\,{th} power, taking
expectations and using Jensen inequality and boundedness of
moments of $Z$, $J_\sigma$ and $T_h J_\sigma$ (the last follows
from the Girsanov theorem, Cauchy--Schwartz inequality and
assumptions on $h$), we get
$$
\ex{\big(T_hZ(t)-Z(t)\big)^{2p}}\le C(p)\left(\int_0^T \abs{h(s)}ds\right)^{2p}.
$$

Further,
\begin{align*}
&\ex{\big(X(t)-\wx(t)\big)^2}\le 3(A_1+A_2+A_3),\\
& A_1 =
\ex{\big(\bj(t)T_{-M\ws\ind{[0,t]}}\big(Z(t)-\wz(t)\big)\big)^2},\\
&A_2=\ex{\big(\big(J_{-\sigma}(t)-\bj(t)\big)T_{-M\ws\ind{[0,t]}}Z(t)\big)^2},\\
&A_3=\ex{\big(J_{-\sigma}(t)\big(
T_{-M\sigma}(1-T_{-M(\ws\ind{[0,t]}-\sigma_t)}\big)Z(t)\big)^2},
\end{align*}
where
\begin{align*}
& J_{-\sigma}(t)=\exp\set{\int_\R M \sigma_t(s) dB^0_s-\frac12 \norm{\sigma_t}^2_{\cH}},\\
& \bj(t)=\exp\set{\int_\R M(\ws\ind{[0,t]})(s)dB^0_s-\frac12
\norm{\ws\ind{[0,t]}}^2_{\cH}}.
\end{align*}
%We first use Cauchy--Shwartz inequality, then Girsanov theorem to
%`cancel' the shift operator, then again Cauchy--Schwartz
%inequality to `isolate' Girsanov density.
Now estimate using Cauchy--Schwartz inequality, Girsanov theorem
(which can be applied as $\si$ and $\ws$ are bounded on $[0,T]$) and
Theorem~\ref{thz1} %!!!!!!!!!!!!!!!!!!!!
\begin{align*}
A_1 & \le
\Big(\ex{\bj^4(t)}\ex{T_{-M\ws\ind{[0,t]}}\big(Z(t)-\wz(t)\big)^4}\Big)^{1/2},\\
& \le C \Big(\ex{\wj(t)\big(Z(t)-\wz(t)\big)^4}\Big)^{1/2}\\
& \le C \Big(\ex{\wj^2(t)}\ex{\big(Z(t)-\wz(t)\big)^8}\Big)^{1/4}\le
C\delta^{2H}.
\end{align*}
Similar reasoning and Lemma~\ref{l1} imply
$$
A_2\le C \ex {\left(\int_\R
M(\ws\ind{[0,t]}-\sigma_t)(s)\,dB^0_s+\frac12\big(\norm{\sigma_t}^2_{\cH}-
\norm{\ws\ind{[0,t]}}^2_{\cH}\big)\right)^2}.
$$
Using condition {\rm (E) 4)}, we obtain $A_2\le C\delta^{2H}$. And
for $A_3$, using the above estimate, we get
$$
A_3\le \int_0^t \abs{M(\ws\ind{[0,t]}-\sigma_t)(s)}ds\le
C\delta^{2H}.
$$
This concludes the proof.
\end{proof}
\begin{remark}\label{simul}
%Note that  it is not always possible to compute $b$ for any $\om$.
%But it is natural to assume that one can compute it for
It is natural to assume that the coefficient $b$ is expressed in
the terms of fBm $B$ rather then in the terms of underlying
Brownian motion $B^0$ (or underlying ``Brownian'' white noise
$\omega$.) This justifies the fact that it is $\sigma$ not
$M\sigma$ what is discretized in \eqref{vocappr}.
\end{remark}
\begin{remark}
Similarly to the proof of Theorem~\ref{th52} one can prove that
for any $s\ge 1$
\begin{equation*}
\ex{\abs{X(t)-\wx(t)}^s}\le \delta^{sH}.
\end{equation*}
The case $s=2$ is considered in the paper to keep classical
``scent'' of results.
\end{remark}
\begin{remark}
Results of this section can be generalized for random initial
condition $X_0$ in the following form: under conditions~{\rm (E)}
and $L^p$-integrability of the initial condition one has convergence
in any $L^{s}$ for $s<p$   with
\begin{equation*}
\ex{\abs{X(t)-\wx(t)}^s}\le \delta^{sH}.
\end{equation*}
Proofs need some simple changes: H\"older inequality for appropriate
powers instead of Cauchy--Schwartz one.
\end{remark}

\end{document}